\documentclass{article}
\usepackage[utf8]{inputenc}
\usepackage{amsmath}
\usepackage{amssymb}
\usepackage{graphicx}
\usepackage{mathtools}
\usepackage[symbol]{footmisc}
\usepackage{array}

\newtheorem{definition}{Definition}[section]
\newtheorem{theorem}{Theorem}[section]
\newtheorem{lemma}[theorem]{Lemma}

\title{ARITHMETIC DIGIT MANIPULATION AND THE CONWAY BASE-13 FUNCTION}
\author{Lyam K. Boylan}
\date{March 2021}

\usepackage[parfill]{parskip}

\begin{document}

\maketitle

\begin{abstract}
    To demonstrate the ability in standard arithmetic operations to perform a variety of digit manipulation tasks, a closed-form representation of the Conway Base-13 Function over the integers is given.
\end{abstract}
\section{Introduction}
Created by the great and late John H. Conway, the Conway Base 13 Function, $f:\mathbb{R}\to\mathbb{R}$, is a counterexample to the converse of the Intermediate Value Theorem. Despite $f$ being discontinuous everywhere, it satisfies that for any interval $(a,b)$, $f$ takes all values between $f(a)$ and $f(b)$. In fact, $f$ takes all values in $\mathbb{R}$ within every interval of non-zero length. 
\begin{figure}[htbp]
\centerline{\includegraphics[scale=.5]{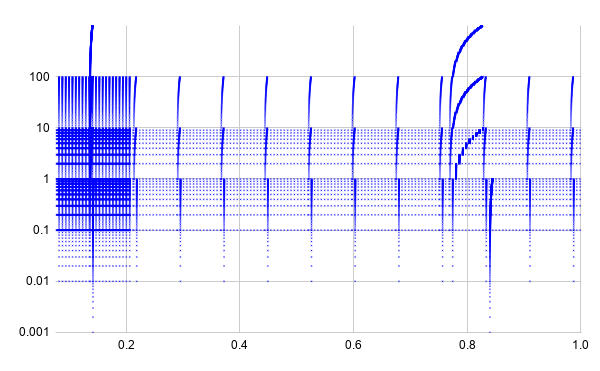}}
\caption{Log-plot of $f$ over a subset of $\mathbb{Z}[\frac{1}{13}]$.}
\label{fig}
\end{figure}

Imperatively, a summary of a definition will be given. Hence, let the set of digits in any base, $b\in\mathbb{Z}_{>1}$, be denoted
$$U_{b}=\{0,\dots,b-1\}.$$
  Thus, the set of decimal and tridecimal digits are $U_{10}=\{0,1,2,3,4,5,6,7,8,9\}$ and $U_{13}=\{0,1,2,3,4,5,6,7,8,9, A,B,C\}$, respectively. The digits $A$, $B$, and $C$ correspond to their decimal equivalents $10$, $11$, and $12$. 
  
 Suppose all $x\in\mathbb{R}^+$ have base-$b$ expansions of the form
    $$x=\dots d_1d_0.d_{-1}d_{-2}\dots_{(b)} \ \text{s.t.} \sum_{k\in \mathbb{Z}}b^kd_k=x$$
where $d_k\in U_{b}$ are individual digits for all $k\in\mathbb{Z}$. Note that $d_k$ corresponds to a digit to the left of the radix point only when $k\leq 0$. In reference to the position of a digit, the term \textit{index} is used. A digit at the $k^{th}$ index of an expansion refers to the digit $k$ positions to the left of the units' column. Hence, a digit at index $0$ is a digit in the units column. If there is no digit at the $k^{th}$ index, the digit is assumed to be zero.
Considering that some values of $x$ and $b$ have two expansions (such as in the cases $0.\overline{9}_{(10)}=1_{(10)}$ or $1.2A\overline{C}_{(13)}=1.2B_{(13)}$), we'll assume the terminating expansion is always preferred. For brevity, we'll introduce the  notation $d_{j\to k(b)}$ as shorthand for $d_jd_{j-1}\dots d_{k+1}d_{k(b)}$. Furthermore, let $d_{j\to k(b)}\subseteq x$  represent, disregarding sign and radix point, that the sequence of digits  $d_{j\to k(b)}$ occurs in the base-$b$ expansion of $x$. E.g;
$$ABC_{(13)}\subseteq -A.BC_{(13)}.$$
If $x\in \mathbb{Z}_{\geq 0}$, then $k< 0\implies d_k=0$. Hence, for integer values, the base-$b$ expansion of $x$ can simply be written $d_{m\to 0}._{(b)}$, where $m\in\mathbb{Z}_{\geq 0}$ is the largest index such that $d_m\neq 0$ (assuming $x\neq 0$, otherwise $m=0$).

Adapted from a definition by Greg Oman \cite{oman}, $f$ is defined in plain language as follows:

For any $x\in\mathbb{R}$, $k\in \mathbb{Z}$, let $d_k$ represent the digit at index $k$ in the tridecimal expansion of $|x|$. A few cases are considered:

\begin{itemize}
\item \textbf{Case 1:} Suppose there exists a digit $A\subset |x|$, such that all digits to the right of such do not contain $A$ or $B$, and there exists exactly one $C\subset |x|$ to the right of such $A$. Let the digits between $A$ and $C$ be denoted $d_{j_A-1\to j_C+1}$, where $j_A$ and $j_C$ are the respective indices of such $A$ and $C$. Let the digits after $C$ be denoted $d_{ j_C-1\to\infty}$. Thus, $f(x)=+d_{j_A-1\to j_C+1}.d_{ j_C-1\to\infty(10)}$.
\item \textbf{Case 2:} Suppose there exists a digit $B\subset|x|$, such that all digits to the right of such do not contain $A$ or $B$, and there exists exactly one $C\subset |x|$ to the right of such $B$. Let the digits between $B$ and $C$ be denoted $d_{j_B-1\to j_C+1}$, where $j_B$ and $j_C$ are the respective indices of such $B$ and $C$. Let the digits after $C$ be denoted $d_{ j_C-1\to\infty}$. Thus, $f(x)=-d_{j_B-1\to j_C+1}.d_{ j_C-1\to\infty(10)}$.
\item \textbf{Otherwise:} $f(x)=0$ if $x$ is not of either form.
\end{itemize}

Here, $d_{j_C-1\to-\infty}$ is shorthand for $\lim_{k\to-\infty}d_{j_C-1\to k}$. It's important to recognize that the final result in cases 1 and 2 are decimal expansions, despite using digits from the tridecimal expansion of $|x|$. This is possible because in either case, the result only uses digits after the right-most $A$ or $B$. Hence, the proceeding digits do not contain $A$ or $B$. There's expectantly exactly one proceeding $C$, (the only other possible tridecimal digit which isn't also a decimal digit) however, which incidentally is excluded in the result. Hence, all digits in the result are indeed decimal. Essentially, $f$ is a recompilation of some of the decimal digits in the tridecimal expansion of $|x|$, using a specific $C$ (if it exists) as a decimal point, and $A$ or $B$ as the sign. Here are a few examples that cover all cases:
\begin{align*}
f(-B1A.3C1415\dots_{(13)})&=\pi\\
f(137_{(13)})&=0\\
f(0.B17C11_{(13)})&=-17.11_{(10)}\\
f(0.\overline{A1C1}_{(13)})&=0\\
f(0.A1\overline{C1}_{(13)})&=0\\
f(0.A999C\overline{9}_{(13)})&=1000_{(10)}
\end{align*}It may be easy to see why $f$ passes through all values of $\mathbb{R}$ within every non-zero-length interval. Regardless, proofs of its properties are not the purpose of this paper. Since $f$ was constructed on the basis of digit manipulation, it lacks a definition using standard mathematical operations.
\begin{theorem}
There exists a closed-form $g:\mathbb{Z}\to\mathbb{R}$ such that $g\subset f$, where $f$ is the Conway Base-13 Function.
\end{theorem}
Understandably, such a prospect would benefit from quantifying its cases. The condition of the existence a digit $A$ or $B\subset |x|$, such that all digits to the right of such do not contain $A$ or $B$ can be quantified as $\exists j_A\big[\ d_{j_A}=A\wedge \not\exists k<j_A(d_{k}\in\{A,B\})\big]$ or $\exists j_A\big[\ d_{j_B}=B\wedge \not\exists k<j_B(d_{k}\in\{A,B\})\big]$ respectively. With the added condition that there exists exactly one $C\subset |x|$ to the right of such $A$ or $B$, the cases become
$$\text{case 1}\iff\exists j_A\Big[\ d_{j_A}=A\wedge \not\exists k<j_A(d_{k}\in\{A,B\})\wedge\exists!j_C<j_A(d_{j_C}=C)\Big]$$
$$\text{case 2}\iff\exists j_B\Big[\ d_{j_B}=B\wedge \not\exists k<j_B(d_{k}\in\{A,B\})\wedge\exists!j_C<j_B(d_{j_C}=C)\Big]$$
This gives rise to an equivalent piecewise formulation:
\[
  f(x) =
  \begin{cases}
                                  + d_{j_A-1\to j_C+1}.d_{j_C-1\to-\infty(10)} & \text{: case 1} \\
                                  - d_{j_B-1\to j_C+1}.d_{j_C-1\to-\infty(10)} & \text{: case 2} \\
  0 & \text{: otherwise}
  \end{cases}
\]
\section{Closed Form Expressions}
As indicated and motivated by Nate Eldredge \cite{wong}, a construction of $f$ using only closed-form arithmetic functions is a possible procedure, albeit tedious and logic-heavy. It requires quite the array of functions designed to arbitrarily manipulate digits and test for logical conditions. 
\subsection{Closed Form Operations}
	As there is no universal definition for closed-form expressions, we assume a conservative definition.
    \begin{definition}
Let an operation be considered closed-form if it can be equivalently expressed in a finite number of operations, of which include addition, subtraction, multiplication, division, exponentiation, principal roots, and the principal branch of the logarithm.
    \end{definition}
     This definition is restrictive so that operations that fulfill this conservative definition expectantly fulfill more liberal ones \cite{chow}. As evident in following sections, the holy-grail of arithmetic digit manipulation relies on the floor and ceiling operations. These can be defined through the use of their relationship to the modulo operation in floored division \cite{knuth}:
 \begin{align*}
\lfloor x\rfloor&:=x-(x\ \text{mod}\  1),\\
\lceil x\rceil&:=x+((-x)\ \text{mod}\ 1).
\end{align*}
Here, $\text{mod}$ is used as a binary operation as opposed to its use in congruence relations. It can be defined though the use of the periodic nature of the principal branch of the logarithm
$$x\ \text{mod}\ y := \frac{y}{2\pi i}\text{Log}\Big(e^{\frac{2\pi i x}{y}}\Big)$$
assuming $0\leq\frac{1}{i}\text{Log}\Big(e^{i\theta}\Big)<2\pi\ \forall\theta\in\mathbb{R}$. Hence, the floor, ceiling, and modulo operations will be considered closed-form. Similarly, the absolute value operation can be defined closed-form through the use of the principal square root, $|x|:=\sqrt{x^2}$. 
\subsection{Logical-Conditional Functions}
Due to the Conway Base-13 Function requiring a piecewise structure, a multitude of  functions that act for testing logical conditions are constructed. In particular: functions that check for equality and inequality relations between two real numbers.
\begin{definition}
Let $E$, “the equivalence function”, be defined as $$E(a,b):=\lfloor(1+\epsilon)^{-|a-b|}\rfloor\ \text{such that}\ \epsilon>0,\ \forall a,b\in\mathbb{R}.$$
 \end{definition} 
It is easily shown that
\[
 E(a,b)=
  \begin{cases}
                                 1 & \text{: $a = b$} \\
                                 0 & \text{: $a \neq b$} \\
  \end{cases}
\]
For brevity, the “negation” of the equivalence function will also be used.
\begin{definition}
Let $N$, “the non-equivalence function”, be defined as $$N(a,b):=1-E(a,b)\ \forall a,b\in\mathbb{R}.$$
 \end{definition}
Similarly,
\[
N(a,b)=
  \begin{cases}
                                 1 & \text{: $a \neq b$} \\
                                 0 & \text{: $a = b$} \\
  \end{cases}
\]
\begin{definition}
Let $G_E$, “the greater-than or equal-to function”, be defined as $$G_E(a,b):=\Big\lfloor\frac{1}{2}+\frac{1}{1+(1+\epsilon)^{b-a}}\Big\rfloor\ \text{such that}\ \epsilon>0,\ \forall a,b\in\mathbb{R}.$$
 \end{definition}
 Although not as trivial as the equivalence function, it can be shown that
 \[
  G_E(a,b)=
  \begin{cases}
                                 1 & \text{: $a \geq b$} \\
                                 0 & \text{: $a <b$} \\
  \end{cases}
\]
 \begin{definition}
Let $M$, “the minimum function”, be defined as $$M(a,b):=aG_E(b,a)+b G_E(a, b)-aE(a, b),\ \forall a,b\in\mathbb{R}.$$
 \end{definition}
With continuing logic, 
\[
 M(a,b)=
  \begin{cases}
                                 a & \text{: $a \leq b$} \\
                                 b & \text{: $a > b$} \\
  \end{cases}
\]
These functions serve as logical-conditional functions, with such, closed-form piecewise-like functions can be crafted arithmetically. 
 \section{Digit Manipulation}
Singling-out digits from an expansion is the most critical ability of digit manipulation. As such, let us introduce the following closed-form functions:
\begin{definition}
Let $\overleftarrow{T}$, “the trailing-digit-truncation function”, be defined as $$\overleftarrow{T}^n_b(x):=\Big\lfloor \frac{x}{b^n}\Big\rfloor$$
$$\forall x\in \mathbb{Z}_{\geq 0},\ \text{for any base}\ b\in \mathbb{Z}_{>1},\ \text{and any digit-index}\ n\in \mathbb{Z}_{\geq0}.$$
\end{definition}
In essence, $\overleftarrow{T}$ removes the right-most $n$ digits from a base-$b$ expansion of $x$. More formally, it removes digits with indices less than a given index $n$.
\begin{lemma}
$x=d_{m\to 0}._{(b)}\implies \overleftarrow{T}_b^n(x)=d_{m\to n}._{(b)}$
\end{lemma}
\textit{Proof:}  Suppose $x=d_{m\to 0}._{(b)}$. By definition of positional notation, $x=\sum_{k=0}^mb^kd_k$. Plugging this into $\overleftarrow{T}$ yields
$$\overleftarrow{T}_b^n(x)=\Bigg\lfloor \frac{\sum_{k=0}^mb^kd_k}{b^n}\Bigg\rfloor=\bigg\lfloor\sum_{k=0}^mb^{k-n}d_k\bigg\rfloor$$
which can be split into a whole and fractional part.
\begin{align*}
&=\bigg\lfloor\sum_{k=n}^mb^{k-n}d_k+\sum_{k=0}^{n-1}b^{k-n}d_k\bigg\rfloor\\
&=\sum_{k=n}^mb^{k-n}d_k+\bigg\lfloor\sum_{k=0}^{n-1}b^{k-n}d_k\bigg\rfloor\\
&=\sum_{k=n}^mb^{k-n}d_k
\end{align*}
We are left with a recompilation of the digits $d_{m\to n}$, such that $d_n$ is now directly to the left of the radix point. In our notation, this is written $d_{m\to n}._{(b)}$.
\hfill $\square$

E.g; $\overleftarrow{T}^2_{10}(123456_{(10)})=1234_{(10)}$. In conjunction, the selection of an arbitrary digit at a given index is possible.

\begin{definition}
Let $D$, “the digit-selection function”, be defined as $$D^n_b(x):=\overleftarrow{T}^n_b(x)-b\overleftarrow{T}^{n+1}_b(x)$$
$$\forall x\in \mathbb{Z}_{\geq 0},\ \text{for any base}\ b\in \mathbb{Z}_{>1},\ \text{and any digit-index}\ n\in \mathbb{Z}_{\geq0}.$$
\end{definition}
This grants the ability to retrieve a digit at the $n^{th}$ index of the base-$b$ expansion of $x$ within a closed-form manner. An ability most critical in construction of the Base-13 Function.
\begin{lemma}
$x=d_{m\to 0}._{(b)}\implies D_b^n(x)=d_n$
\end{lemma}
\textit{Proof:}  Suppose $x=d_{m\to 0}._{(b)}$. Using the result in our previous lemma, $D$ becomes
\begin{align*}
D^n_b(x)&=\sum_{k=n}^mb^{k-n}d_k-b\sum_{k=n+1}^mb^{k-n-1}d_k\\
&=\sum_{k=n}^mb^{k-n}d_k-\sum_{k=n+1}^mb^{k-n}d_k\\
\end{align*}
which cleanly cancels all digits with indices other than $k=n$
\begin{align*}
&=d_n+\sum_{k=n+1}^mb^{k-n}d_k-\sum_{k=n+1}^mb^{k-n}d_k\\
&=d_n
\end{align*}
Leaving our desired digit, $d_n$.
\hfill $\square$

E.g; $D^2_{10}(123456_{(10)})=4_{(10)}$. Not surprisingly, the amount of digits in an expansion can also be deduced arithmetically.

\begin{definition}
Let L, “the length function”, be defined as 
$$L_b(x):=\lceil \log_b(x+1)\rceil$$
$$\forall x\in \mathbb{Z}_{\geq 0},\ \text{for any base}\ b\in \mathbb{Z}_{>1}.$$
\end{definition}
This is variant to the usual method to count the number of digits: $\lfloor\log_b(x)\rfloor+1$. However the latter is undefined for the case $x=0$, whereas $L_b(0)=0$. Otherwise both methods are equivalent over the positive integers.
\begin{lemma}
$x=d_{m\to 0}._{(b)}\wedge x>0\implies L_b(x)=m+1$
\end{lemma}
\textit{Proof}:  Suppose $x=d_{m\to 0}.{(b)}\wedge x>0$,
\begin{align*}
\implies L_b(x)&=\bigg\lceil \log_b\Big(1+\sum_{k=0}^mb^k d_k\Big)\bigg\rceil\\
\implies \Big\lceil \log_b(b^m)\Big\rceil<L_b(x)&\leq\Big\lceil \log_b(b^{m+1})\Big\rceil\\
\implies m<L_b(x)&\leq m+1\\
\implies L_b(x)&=m+1
\end{align*}
\hfill $\square$

E.g; $L_{10}(10_{(10)})=L_{10}(99_{(10)})=2$. Using  the length function, $L$, to bound a summation of the equivalence function, $E$, in conjunction with the digit-selection function, $D$, grants a quaint method to count the number of occurrences of a digit in a given number.

\begin{definition}
Let $O$, “the digit-occurrence-counting function”, be defined as 
$$O_b^p(x):=\sum_{k=0}^{L_b(x)-1}E(D_b^k(x),p)$$
$$\forall x\in \mathbb{Z}_{\geq 0},\ \text{for any base}\ b\in \mathbb{Z}_{>1},\ \text{and any digit}\ p\in U_b.$$
\end{definition}

It should be evident that as $O$ loops through all possible indices, $k$, for digits in the base-$b$ expansion of $x$, the summation increments by 1 iff the digit at index $k$ is equivalent to the given digit $p$, which we are looking to count the occurrences of. In other words, $O$ counts the number of occurrences of a digit $p$ in the base-$b$ expansion of $x$.

Less trivial is a method to deduce the a specific index of an occurrence of a given digit.

\begin{definition}
Let I, “the digit-occurrence-index function”, be defined as 
$$I_b^p(x):=\sum_{k=1}^{L_b(x)}E\Big(O_b^p\big(\overleftarrow{T}_b^k(x)\big),O_b^p(x)\Big)$$
$$\forall x\in \mathbb{Z}_{\geq 0},\ \text{for any base}\ b\in \mathbb{Z}_{>1},\ \text{and any digit}\ p\in U_b.$$
\end{definition}
The purpose of $I$ is to return the index of the right-most digit $p$ in the base-$b$ expansion of $x$. If there isn't such an index, then $I$ returns $L_b(x)$, which is by definition, a number higher than the maximum index of a nonzero digit.
\begin{lemma}
\[
x=d_{m\to 0}._{(b)}\implies I_b^p(x)=
  \begin{cases}
                                j  & \textup{: $\exists j\big[d_j =p\wedge\forall k<j(d_k\neq p)\big]$} \\
                                L_b(x) & \textup{: otherwise} \\
  \end{cases}
\]
\end{lemma}
\textit{Proof:} Suppose $x=d_{m\to 0}._{(b)}$. We'll look at the case where there does exist a right-most digit $p$ in the base-$b$ expansion of $x$. 

\textbf{Case 1: }  $\exists j\big[d_j =p\wedge\forall k<j(d_k\neq p)]$

Thus, with such a digit having index $j$, truncating off digits of $x$ with indices less than $k$ for $k\leq j$, yields a number with no occurrences of $p$ removed. Likewise, truncating for $k>j$ yields a number with at least one less occurrence of $p$.
\begin{align*}
k\leq j\iff O_b^p\big(\overleftarrow{T}_b^k(x)\big)&=O_b^p(x)\\
\implies E\Big(O_b^p\big(\overleftarrow{T}_b^k(x)\big),O_b^p(x)\Big)&=
  \begin{cases}
                                1  & \text{: $k \leq j$} \\
                                0 & \text{: $k >j$}
  \end{cases}
\end{align*}
Thus, the summation can be split into
$$I_b^p(x)=\sum_{k=1}^{j}1+\sum_{k=j+1}^{L_b(x)}0=j$$
Resulting in the index, $j$.

\textbf{Case 2:} $\nexists j\big[d_j =p\wedge\forall k<j(d_k\neq p)\big]$

In the other case, since $x$ is an integer of finite digits, there not being a right-most digit $p$ implies that there are no occurrences.
\begin{align*}
 O_b^p\big(\overleftarrow{T}_b^k(x)\big)=O_b^p(x)&=0\ \forall k\\
\implies E\Big(O_b^p\big(\overleftarrow{T}_b^k(x)\big),O_b^p(x)\Big)&=1\\
\implies I_b^p(x)&=\sum_{k=1}^{L_b(x)}1\\
 &=L_b(x)
\end{align*}As such, the sum is trivially the bound, $L_b(x)$.
\[
\therefore x=d_{m\to 0}._{(b)}\implies I_b^p(x)=
  \begin{cases}
                                j  & \text{: $\exists j\big[d_j =p\wedge\forall k<j(d_k\neq p)\big]$} \\
                                L_b(x) & \text{: otherwise} \\
  \end{cases}
\]
\hfill $\square$

In parody to $\overleftarrow{T}$, we'll define a function that virtually removes the left-most $n$ digits from a base-$b$ expansion of $x$.
\begin{definition}
Let $\overrightarrow{T}$, “the leading-digit-truncation function”, be defined as $$\overrightarrow{T}^n_b(x):=\sum_{k=0}^{L_b(x)-n-1}b^kD_b^k(x)$$
$$\forall x\in \mathbb{Z}_{\geq 0},\ \text{for any base}\ b\in \mathbb{Z}_{>1},\ \text{for any digit-index}\ n\in \mathbb{Z}_{\geq 0}$$
\end{definition}
Clearly, $\overrightarrow{T}$ reassembles the digits in the base-$b$ expansion of $x$ into their original position, save for the last $n$ digits.
\begin{definition}
Let $K$, the “cut-to-index function” be defined as 
$$K_b^p(x):=\sum_{k=0}^{I_b^p(x)}b^kD_b^k(x)$$
$$\forall x\in \mathbb{Z}_{\geq 0},\ \text{for any base}\ b\in \mathbb{Z}_{>1},\ \text{for any digit}\ p\in U_b.$$
\end{definition}
Similar to $\overrightarrow{T}$, $K$ reassembles the digits in the base-$b$ expansion of $x$ into their original position, save for the last digits with indices greater than $I_b^p(x)$. For the case where $p\nsubseteq x$, we find that $I_b^p(x)=L_b(x)$, which implies that $K_b^d(x)=x$.

\section{Assembling The Conway Base-13 Function}

Perhaps the most daunting of tasks to replicate in the Conway Base-13 Function is recompiling digits in an expansion from one base to another, and replacing a digit with a radix-point.

\begin{definition}
Let $X$, “the base-to-base re-radix function”, be defined as 
$$X_{b_1, b_2}^p(x):=\sum_{k=0}^{L_{b_1}(x)-1}N(D_{b_1}^k(x), p)D_{b_1}^k(x)b_2^{k-I_{b_1}^p(x)-G_E(I_{b_1}^p(x), k)}$$
$$\forall x\in \mathbb{Z}_{\geq 0},\ \text{for any bases}\ b_1, b_2\in \mathbb{Z}_{>1},\ \text{and any digit}\ p\in U_p.$$
\end{definition}
$X$ removes a specific digit $p$, with index $j$ in the base-$b_1$ expansion of $x$. This position will be virtually used as a new radix-point. Digits to the left of $p$ (with indices $k>j$) are placed directly to left of this new radix, and digits to the right of $p$ (with indices $k<j$) are placed directly to the right. The final result is treated as a base-$b_2$ expansion. Although still returning defined values for the instances where there are multiple occurrences of $p$, its behaviour is irrelevant, as such a case is evidently disregarded in further construction of the Base-13 function.
\begin{lemma}
$$x=d_{m\to 0}._{(b_1)}\wedge\exists!j(d_j=p)\implies X_{b_1, b_2}^p(x)=d_{m\to j+1}.d_{j-1\to 0(b_2)}\ \forall b_2\in\mathbb{Z}_{\geq b_1}$$
\end{lemma}
\textit{Proof}: Suppose $x=d_{m\to 0}._{(b_1)}$ and $\exists!j(d_j=p)$. Thus the index, $j$, is given by $I_{b_1}^p(x)=j$. A digit at index $k$ is given by $D_{b_1}^k(x)=d_k$. Hence $\forall b_2\in\mathbb{Z}_{\geq b_1}$, substituting for our positional notation,
\[
N(D_{b_1}^k(x), p)D_{b_1}^k(x)b_2^{k-I_{b_1}^p(x)-G_E(I_{b_1}^p(x), k)}=
  \begin{cases}
                                  d_k b_2^{k-j} & \text{if $ k < j$} \\
                                  0 & \text{if $ k = j$} \\
                                  d_k b_2^{k-j-1} & \text{if $ k > j$} \\
  \end{cases}\\
\]
which can be used to split the sum into
$$ X_{b_1, b_2}^p(x) = \sum_{k=0}^{j-1}d_kb_2^{k-j}+\sum_{k=j+1}^{m}d_kb_2^{k-j-1}$$
We are left with two recompilations of digits from base-$b_1$ to base-$b_2$, with the digits to the left of $p$ directly to left of the radix, and digits to the right of $p$ to the right. In our positional notation, this is equivalent to $d_{m\to j+1}.d_{j-1\to 0(b_2)}$.

\hfill $\square$

E.g; $X_{13, 10}^{C}(1C3_{(13)})=1.3_{(10)}$. Lastly, we'll introduce a method to detect whether one of two given digits are contained within a base-$b$ expansion. This will act as the step in determining if the final expansion of Conway’s Base-13 function will be positive or negative.

\begin{definition}
Let $S$, “the resulting-sign function”, be defined as 

$$S_b^{p_1,p_2}(x):=E\big(O_b^{p_1}(x),1\big)-E\big(O_b^{p_2}(x),1\big)$$ 
$$\forall x\in \mathbb{Z}_{\geq 0},\ \text{for any base}\ b\in \mathbb{Z}_{>1},\ \text{for any digits}\ p_1, p_2\in \{0,\dots,b-1\}.$$
\end{definition}
Unlike the previous function, $S$ is much simpler in description. If there exists exactly one $p_1\subset x$, and not exactly one $p_2\subset x$ (assuming a base-$b$ expansion), then $S_b^{p_1,p_2}(x)=1$. Similarly, $S_b^{p_1,p_2}(x)=-1$ if there exists exactly one $p_2\subset x$, and not exactly one $p_1\subset x$. Otherwise the result is zero.
\begin{lemma}
\[
x=d_{m\to 0}._{(b)}\implies S_b^{p_1, p_2}(x)=
  \begin{cases}
                                +1  & \textup{: $\exists! j_1(d_{j_1} =p_1)\wedge\nexists! j_2(d_{j_2} =p_2)$} \\
                                -1  & \textup{: $\nexists! j_1(d_{j_1} =p_1)\wedge\exists! j_2(d_{j_2} =p_2)$} \\
                                0  & \textup{: otherwise} \\
  \end{cases}
\]
\end{lemma}
\textit{Proof}: Suppose $x=d_{m\to 0}._{(b)}$. With the definitions of $E$ and $O$, the values of $S$, defined by $E\big(O_b^{p_1}(x),1\big)-E\big(O_b^{p_2}(x),1\big)$, in the following case-table are straightforward.
\begin{center}
\begin{tabular}{ |c||c|c| } 
 \hline
cases & $\exists! j_1(d_{j_1} =p_1)$&$\nexists! j_1(d_{j_1} =p_1)$ \\ 
 \hline
 \hline
$\exists! j_2(d_{j_2} =p_1)$ & $S_b^{p_1,p_2}(x)=\ 1-1\ =\ 0$ & $S_b^{p_1,p_2}(x)=\ 0-1\ =\ -1$ \\ 
\hline
$\nexists! j_2(d_{j_2} =p_2)$ & $S_b^{p_1,p_2}(x)=\ 1-0\ =\ 1$ & $S_b^{p_1,p_2}(x)=\ 0-0\ =\ 0$ \\ 
 \hline
\end{tabular}
\end{center}
\hfill $\square$

With an arsenal of closed-form logical-conditional and digit manipulating functions, the Conway Base-13 Function can be constructed in three phases.
\subsection{Phase 1}

$$f_1(x)=M\Big(K_{13}^{A}|x|, K_{13}^{B}|x|\Big)$$
In the initial phase, all digits directly to the left of the rightmost $A$ or $B$ in the tridecimal expansion of $x$ are truncated. As the sign of the input is disregarded in Conway's Base-13 Function, $f$, the absolute value of $x$ is taken for each reference of $x$ in $f_3$. For any $k\in\mathbb{Z}$, let $d_k$ represent the digit at index $k$ in the tridecimal expansion of $|x|$. Let the rightmost-index of $A$ be written as $I_{13}^{A}|x|=j_A$ and the rightmost-index of $B$ be written as $I_{13}^{B}|x|=j_B$. Note that by our definitions of $I$ and $K$,
$$A\nsubseteq x\implies j_A=L_{13}|x|\implies K_{13}^{A}|x|=|x|$$
$$B\nsubseteq x\implies j_B=L_{13}|x|\implies K_{13}^{B}|x|=|x|$$
Since $L$ is monotonically increasing, the inequality relation between $j_A$ and $j_B$ implies
$$j_A \leq j_B \iff L_{13}\Big(K_{13}^{A}|x|\Big) \leq L_{13}\Big(K_{13}^{B}|x|\Big)\iff K_{13}^{A}|x|\leq K_{13}^{B}|x|$$
$$j_B \leq j_A \iff L_{13}\Big(K_{13}^{B}|x|\Big) \leq L_{13}\Big(K_{13}^{A}|x|\Big)\iff K_{13}^{B}|x|\leq K_{13}^{A}|x|$$
Therefore, they determine the value of $M$ by
$$j_A\leq j_B \iff M\Big(K_{13}^{A}|x|, K_{13}^{B}|x|\Big)=K_{13}^{A}|x|$$
$$j_B\leq j_A \iff M\Big(K_{13}^{A}|x|, K_{13}^{B}|x|\Big)=K_{13}^{B}|x|$$
Consequently $f_1$ becomes
\[
  f_1(x) =
  \begin{cases}
                                  K_{13}^{A}|x| & \text{: $j_A<j_B$} \\
                                  K_{13}^{B}|x| & \text{: $j_B<j_A$} \\
 								  |x| & \text{: otherwise} \\
  \end{cases}
\]
And by definition of $K$, can be represented
\[
  f_1(x) =
  \begin{cases}
                                  Ad_{j_A-1\to 0}._{(13)} & \text{: $A\subseteq x\wedge j_A<j_B$} \\
                                  Bd_{j_B-1\to 0}._{(13)} & \text{: $B\subseteq x\wedge j_B<j_A$} \\
 								  |x| & \text{: otherwise} \\
  \end{cases}
\]
\subsection{Phase 2}
$$f_2(x)=f_1(x)E\Big(O_{13}^{C}\big(f_1(x)\big), 1\Big)$$
In phase 2, we are checking if after such an $A$ or $B$, there exists exactly one $C$ leftover from phase 1. If there does not exist exactly one such $C$,
$$O_{13}^{C}\big(f_1(x)\big)\neq 1 \implies E\Big(O_{13}^{C}\big(f_1(x)\big), 1\Big) = 0 \implies f_2(x)=0$$
Otherwise, let the index of such be denoted $I_{13}^{C}\big(f_1(x)\big)=j_C$. Therefore
\[
  f_2(x) =
  \begin{cases}
                                  Ad_{j_A-1\to j_C+1} C d_{j_C-1\to 0}._{(13)} & \text{: $A\subseteq f_1(x)\wedge C\subseteq f_1(x)$} \\
                                   Bd_{j_B-1\to j_C+1} C d_{j_C-1\to 0}._{(13)} & \text{: $B\subseteq f_1(x)\wedge C\subseteq f_1(x)$} \\
 								   0 & \text{: $C\nsubseteq f_1(x)$} \\
                                   f_1(x) & \text{: otherwise} \\
  \end{cases}
\]
\subsection{Phase 3}
$$f_3(x)=S_{13}^{A, B}\Big(f_2(x)\Big)X_{13, 10}^{C}\Big(\overrightarrow{T}^1_{13}\big(f_2(x)\big)\Big)$$
The last phase determines the sign of the result, truncates off the leftover $A$ or $B$, recompiles the tridecimal expansion into decimal, and essentially replaces $C$ with a decimal point (assuming $C\subseteq f_2(x)$). By the definition of $S$ and $f_2$, we find that
\[
S_{13}^{A, B}\Big(f_2(x)\Big) =
  \begin{cases}
                                  +1 & \text{: $A\subseteq f_2(x)$} \\
                                  -1 & \text{: $B\subseteq f_2(x) $} \\
  0 & \text{: otherwise}
  \end{cases}
\]
and by our definition of $\overrightarrow{T}$,
\[
 \overrightarrow{T}^1_{13}\big(f_2(x)\big) =
  \begin{cases}
                                  d_{j_A-1\to j_C+1} C d_{j_C-1\to 0}._{(13)} & \text{: $A\subseteq f_2(x)\wedge C\subseteq f_2(x)$} \\
                                    d_{j_B-1\to j_C+1} C d_{j_C-1\to 0}._{(13)} & \text{: $B\subseteq f_2(x)\wedge C\subseteq f_2(x)$} \\
 								   0 & \text{:  $C\nsubseteq f_2(x)$} \\
 								   d_{L_{13}(f_2(x))-2\to j_C+1}Cd_{j_C-1\to 0}._{(13)} & \text{:  otherwise} \
  \end{cases}
\]
Hence, by definition of $X$\footnote[1]{Note that the value in the final case is irrelevant, as such a case results in 0 in $f_3$.},
\begin{gather*} 
X_{13, 10}^{C}\Big(\overrightarrow{T}^1_{13}\big(f_2(x)\big)\Big) = \\
 \begin{cases}
                                  d_{j_A-1\to j_C+1} . d_{j_C-1\to 0(10)} & \text{: $A\subseteq f_2(x)\wedge C\subseteq f_2(x)$} \\
                                    d_{j_B-1\to j_C+1} . d_{j_C-1\to 0(10)} & \text{: $B\subseteq f_2(x)\wedge C\subseteq f_2(x)$} \\
 								   0 & \text{: $C\nsubseteq f_2(x)$} \\
                                    d_{L_{13}(f_2(x))-2\to j_C+1}.d_{j_C-1\to 0(10)} & \text{: otherwise} \\
  \end{cases}
\end{gather*}
Therefore, the final result is of form
\[
  f_3(x) =
  \begin{cases}
                                  +d_{j_A-1\to j_C+1} . d_{j_C-1\to 0}._{(10)} & \text{: $A\subseteq f_2(x)\wedge C\subseteq f_2(x)$} \\
                                   -d_{j_B-1\to j_C+1} . d_{j_C-1\to 0}._{(10)} & \text{: $B\subseteq f_2(x)\wedge C\subseteq f_2(x)$} \\
 								   0 & \text{: otherwise} \\
  \end{cases}
\]
It should be seen that the values for each case in $f_3$ is equivalent to $f$ when $x$ is an integer. Furthermore, as the cases from the original quantification hold equivalently, we find that $f_3\subset f$, directly satisfying that $f_3$ is a closed-form representation of the Conway Base-13 Function over the integers.
$$\text{QED}$$
\section{Concluding Remarks}
Due to the fractal-like symmetry of $f$, such that $f(13^nx)=f(x) \forall n\in\mathbb{Z}$, it's possible to extend $f_3:\mathbb{Z}\to\mathbb{R}$ to $f_3:\mathbb{Z}[\frac{1}{13}]\to\mathbb{R}$ by imposing that if $x=\frac{y}{13^n} \forall y\in\mathbb{Z},\ \forall n\in\mathbb{Z}_{\geq 0}$, then $f_3(x)=f_3(y)$; as was done in the creation of figure \ref{fig}. It may be possible to extend $f_3$ to even larger sets of numbers who's distribution of digits are known, however a closed-form for $f$ over the entirety of $\mathbb{R}$ is impossible, as the digit-distrubution for every real number is not computable \cite{wong}.

It is no doubt that the computational efficiency of these algorithms are far from optimal. Of course, a computer can perform a variety of digit manipulation tasks directly and quite efficiently, without the use of arithmetic closed-form functions such as these. The purpose, rather, was to fulfill the recreational endeavour of finding the first equation for Conway's Base-13 Function, based solely on finite arithmetic; inspired to completion by the recent passing of John H. Conway (1937$\to$2020).

\textbf{Acknowledgement} Humblest gratitude to my mother, Loral Boylan, and to Dr.  Peter Dukes for suggestions on improving the readability of this paper.

 \end{document}